\documentclass[12pt]{article}
\usepackage{amsthm,amsfonts}
\usepackage{fullpage}

\newtheorem{theorem}{Theorem}
\newtheorem{corollary}[theorem]{Corollary}
\newtheorem{lemma}[theorem]{Lemma}

\title{Words avoiding $\frac{7}{3}$-powers and the Thue-Morse morphism}
\author{Narad Rampersad \\
School of Computer Science \\
University of Waterloo \\
Waterloo, ON, N2L 3G1 \\
CANADA \\
{\tt nrampersad@math.uwaterloo.ca}}

\begin{document}
\date{\today}
\maketitle

\begin{abstract}
In 1982, S\'{e}\'{e}bold showed that the only overlap-free binary
words that are the fixed points of non-identity morphisms are the
Thue-Morse word and its complement.  We strengthen S\'{e}\'{e}bold's
result by showing that the same result
holds if the term `overlap-free' is replaced with
`$\frac{7}{3}$-power-free'.  Furthermore, the number $\frac{7}{3}$
is best possible.
\end{abstract}

\section{Introduction}
Let $\Sigma$ be a finite, non-empty set called an \emph{alphabet}.
We denote the set of all finite words over the alphabet
$\Sigma$ by $\Sigma^*$.  We also write $\Sigma^+$ to denote
the set $\Sigma^*-\{\epsilon\}$, where $\epsilon$ is the empty word.
Let $\Sigma_k$ denote the alphabet $\{0,1,\ldots,k-1\}$.
Throughout this paper we will work exclusively with the binary alphabet
$\Sigma_2$.

Let $\mathbb{N}$ denote the set $\{0,1,2,\ldots\}$.
An \emph{infinite word} is a map from
$\mathbb{N}$ to $\Sigma$, and a \emph{bi-infinite} word is a map from
$\mathbb{Z}$ to $\Sigma$.  The set of all infinite words over the
alphabet $\Sigma$ is denoted $\Sigma^\omega$.  We also write
$\Sigma^\infty$ to denote the set $\Sigma^* \cup \Sigma^\omega$.

A map $h:\Sigma^* \rightarrow \Delta^*$ is called a \emph{morphism} if
$h$ satisfies $h(xy)=h(x)h(y)$ for all $x,y\in\Sigma^*$.  A morphism
may be defined simply by specifying its action on $\Sigma$.
A morphism $h:\Sigma^* \rightarrow \Sigma^*$ such that $h(a)=ax$
for some $a\in\Sigma$ is said to be \emph{prolongable on $a$}; we
may then repeatedly iterate $h$ to obtain the \emph{fixed point}
$h^\omega(a)=axh(x)h^2(x)h^3(x)\cdots$.

An \emph{overlap} is a word of the form $axaxa$, where $a\in\Sigma$
and $x\in\Sigma^*$.  A word $w'$ is called a subword of $w\in\Sigma^\infty$
if there exist $u\in\Sigma^*$ and $v\in\Sigma^\infty$ such that $w=uw'v$.
We say a word $w$ is \emph{overlap-free} (or \emph{avoids overlaps})
if no subword of $w$ is an overlap.

Let $\mu$ be the \emph{Thue-Morse morphism}; \emph{i.e.}, the morphism
defined by $\mu(0)=01$ and $\mu(1)=10$.  It is well-known \cite{MH44, Thu12}
that the \emph{Thue-Morse word}, $\mu^\omega(0)$, is overlap-free.
The properties of overlap-free words have been studied extensively
(see, for example, the survey by S\'{e}\'{e}bold \cite{See84}).
S\'{e}\'{e}bold \cite{See82,See85} showed that $\mu^\omega(0)$ and
$\mu^\omega(1)$ are the only infinite overlap-free binary words that can be
obtained by iteration of a morphism.  Another proof of this fact was later
given by Berstel and S\'{e}\'{e}bold \cite{BS93}.
We will show that this result can be strengthened somewhat.
We will first need the notion of a \emph{fractional power},
which was first introduced by Dejean \cite{Dej72}.

Let $\alpha$ be a rational number such that $\alpha \geq 1$.
An \emph{$\alpha$-power} is a word of the form $x^nx'$, where
$x,x'\in\Sigma^*$, and $x'$ is a prefix of $x$ with $n+|x'|/|x|=\alpha$.
We say a word $w$ is \emph{$\alpha$-power-free} (or
\emph{avoids $\alpha$-powers}) if no subword of $w$ is an
$\beta$-power for any rational $\beta\geq\alpha$; otherwise, we say
$w$ \emph{contains an $\alpha$-power}.  Note that a word is overlap-free
if and only if it is $(2+\epsilon)$-power-free for all $\epsilon > 0$;
for example, an overlap-free word is necessarily $\frac{7}{3}$-power-free.

In this paper we will be particularly concerned with $\frac{7}{3}$-powers.
Several results previously known for overlap-free binary words have recently
been shown to be true for $\frac{7}{3}$-power-free binary words as well.
For example, Restivo and Salemi's factorization theorem for overlap-free
binary words \cite{RS84} was recently shown to be true for
$\frac{7}{3}$-power-free binary words by Karhum\"{a}ki and Shallit
\cite{KS03}.  In 1964, Gottschalk and Hedlund \cite{GH64} showed that
the bi-infinite overlap-free binary words were simply shifts of the bi-infinite
analogue of the Thue-Morse word, and in 2000, Shur \cite{Shu00} showed
that a similar result holds for the bi-infinite $\frac{7}{3}$-power-free
binary words.  Furthermore, Shur showed that the number $\frac{7}{3}$ is
best possible.

The goal of this paper is to generalize S\'{e}\'{e}bold's result by
showing that $\mu^\omega(0)$ and $\mu^\omega(1)$ are the only infinite
$\frac{7}{3}$-power-free binary words that can be obtained by iteration
of a morphism.  At first glance, it may seem that this is an immediate
consequence of Shur's result; however, this is not necessarily so, as
there are infinite $\frac{7}{3}$-power-free binary words that cannot be
extended to the left to form bi-infinite $\frac{7}{3}$-power-free
binary words.  For example, the infinite binary word $001001\mu^\omega(1)$,
which was shown by Allouche \emph{et al.} \cite{ACS98} to be the
lexicographically least infinite overlap-free binary word, cannot
be extended to the left to form a $\frac{7}{3}$-power-free
word:  prepending a 0 creates the cube 000, and prepending a 1 creates
the $\frac{7}{3}$-power 1001001.

\section{Preliminary lemmata}
We will need the following result due to Shur \cite{Shu00}.

\begin{theorem}[Shur]
\label{KS1}
Let $w\in\Sigma_2^*$, and let $\alpha > 2$ be a real number.  Then $w$
is $\alpha$-power-free iff $\mu(w)$ is $\alpha$-power-free.
\end{theorem}

We will also make frequent use of the following result due to
Karhum\"{a}ki and Shallit \cite{KS03}.  This theorem is a
generalization of a similar factorization theorem for overlap-free words
due to Restivo and Salemi \cite{RS84}.

\begin{theorem}[Karhum\"{a}ki and Shallit]
\label{KS2}
Let $x\in\Sigma_2^*$ be a word avoiding $\alpha$-powers, with
$2 < \alpha \leq \frac{7}{3}$.  Then there exist $u,v,y$ with
$u,v\in\{\epsilon,0,1,00,11\}$ and a word $y\in\Sigma_2^*$
avoiding $\alpha$-powers, such that $x=u\mu(y)v$.
\end{theorem}

Next, we will establish a few lemmata.  Lemma~\ref{subw} is analogous
to a similar lemma for overlap-free words given in Allouche and Shallit
\cite[Lemma 1.7.6]{AS03}.  (This result was also stated without
formal proof by Berstel and S\'{e}\'{e}bold \cite{BS93}.)

\begin{lemma}
\label{subw}
Let $w\in\Sigma_2^*$ be a $\frac{7}{3}$-power-free word with $|w| \geq 52$.
Then $w$ contains $\mu^3(0)=01101001$ and $\mu^3(1)=10010110$ as subwords.
\end{lemma}

\begin{proof}
Since $w$ is $\frac{7}{3}$-power-free, by Theorem~\ref{KS2} we can
write
\begin{equation}
w=u\mu(y)v,
\end{equation}
where $y$ is $\frac{7}{3}$-power-free and $|y| \geq 24$.  Similarly,
we can write
\begin{equation}
y=u'\mu(y')v',
\end{equation}
where $y'$ is $\frac{7}{3}$-power-free and $|y'| \geq 10$.  Again,
we can write
\begin{equation}
y'=u''\mu(y'')v'',
\end{equation}
where $y''$ is $\frac{7}{3}$-power-free and $|y''| \geq 3$.  From Equations
(1)--(3), we get
\begin{eqnarray*}
w & = & u\mu(u'\mu(u''\mu(y'')v'')v')v \\
  & = & u\mu(u')\mu^2(u'')\mu^3(y'')\mu^2(v'')\mu(v')v,
\end{eqnarray*}
where $u,u',u'',v,v',v''\in\{\epsilon,0,1,00,11\}$.  Since
$y''$ is $\frac{7}{3}$-power-free and $|y''| \geq 3$, $y''$ contains
both 0 and 1, and so $\mu^3(y'')$, and consequently $w$, contains both
$\mu^3(0)=01101001$ and $\mu^3(1)=10010110$ as subwords as required.
\end{proof}

\begin{lemma}
\label{abb}
Let $w'$ be a subword of $w\in\Sigma_2^*$, where $w'$ is either
of the form $abb\mu(w'')$ or $\mu(w'')bba$ for some $a,b\in\Sigma_2$
and $w''\in\Sigma_2^*$.  Suppose also that $a \neq b$ and
$|w''| \geq 2$.  Then $w$ contains a $\frac{7}{3}$-power.
\end{lemma}

\begin{proof}
Suppose $ab=10$ and $w'=100\mu(w'')$ (the other cases follow similarly).
The word $\mu(w'')$ may not begin with a 0 as that would create the
cube 000.  Hence we have $w'=10010\mu(w''')$ for some
$w'''\in\Sigma_2^*$.  If $\mu(w''')$ begins with 01, then
$w'$ contains the $\frac{7}{3}$-power 1001001.  If $\mu(w''')$
begins with 10, then $w'$ contains the $\frac{5}{2}$-power 01010.
Hence, $w$ contains a $\frac{7}{3}$-power.
\end{proof}

\begin{lemma}
\label{pattern}
For $i,j\in\mathbb{N}$, let $w$ be a $\frac{7}{3}$-power-free
word over $\Sigma_2$ such that $|w| = (7+2j)2^i - 1$.
Let $a$ be an element of $\Sigma_2$.  Then $waw$ contains a
$\frac{7}{3}$-power $x$, where $|x| \leq 7 \cdot 2^i$.
\end{lemma}

\begin{proof}
Suppose $a=1$ (the case $a=0$ follows similarly).
The proof is by induction on $i$.  For the base case we have
$i=0$.  Hence, $|w| \geq 6$ and $|w|$ is even.
If $w$ either begins or ends with 11, then $w1w$ contains
the cube 111, and the result follows.  Suppose then that $w$ neither begins
nor ends with 11.  By explicitly examining all 13 words of length
six that avoid $\frac{7}{3}$-powers and neither begin nor end with 11,
we see that all such words of length at least six can be written in the
form $pbbq$, where $p,q\in\Sigma_2^+$ and $b\in\Sigma_2$.  Hence, $w1w$ must
have at least one subword with prefix $bb$ and suffix $bb$.  Moreover,
since $|w|$ is even, there must exist such a subword where the prefix $bb$
and the suffix $bb$ each begin at positions of different parity in $w1w$.
Let $x$ be a smallest such subword such that $w1w$ neither
begins nor ends with $x$.  Suppose $b=0$ (the case $b=1$ follows similarly).
Then $x=000$, $x=00100$, or $x$ contains a subword $01010$ or $10101$.
Hence, $w1w$ contains one of the subwords 000, 01010, 10101, or 1001001
as required.

Let us now assume that the lemma holds for all $i'$, where $0<i'<i$.
Since $w$ avoids $\frac{7}{3}$-powers, and since $|w| \geq 7$, by
Theorem~\ref{KS2} we can write $w=u\mu(w')v$, where
$u,v\in\{\epsilon,0,1,00,11\}$ and $w'\in\Sigma_2^*$ is
$\frac{7}{3}$-power-free.  By applying
a case analysis similar to that used in Cases (1)--(4) of the
proof of Theorem~\ref{h_eq_mu} below, we can eliminate all but
three cases: $(u,v)\in\{(\epsilon,\epsilon),(\epsilon,0),(0,\epsilon)\}$.

\begin{enumerate}
\renewcommand{\labelenumi}{Case \arabic{enumi}:}
\item
$(u,v)=(\epsilon,\epsilon)$.
In this case $w=\mu(w')$.  This is clearly not possible, since for
$i>0$, $|w| = (7+2j)2^i - 1$ is odd.

\item
$(u,v)=(\epsilon,0)$.
Then $w=\mu(w')0$ and $w1w=\mu(w')01\mu(w')0=\mu(w'0w')0$.  If
$|w| = (7+2j)2^i - 1$, we see that
$|w'| = (7+2j)2^{i-1} - 1$.  Hence, if $i'=i-1$, we may
apply the inductive assumption to $w'0w'$.  We thus obtain that
$w'0w'$ contains a $\frac{7}{3}$-power $x'$, where $|x'| \leq 7 \cdot 2^{i-1}$,
and so $w1w$ must contain a $\frac{7}{3}$-power $x=\mu(x')$, where
$|x| \leq 7 \cdot 2^i$.

\item
$(u,v)=(0,\epsilon)$.  This case is handled similarly to the previous
case, and we omit the details.
\end{enumerate}

By induction then, we have that $waw$ contains a
$\frac{7}{3}$-power $x$, where $|x| \leq 7 \cdot 2^i$.
\end{proof}

\begin{lemma}
\label{pattern2}
For $i\in\mathbb{N}$, let $w$ be a $\frac{7}{3}$-power-free
word over $\Sigma_2$ such that $|w| = 5 \cdot 2^i - 1$.
Let $a$ be an element of $\Sigma_2$.  Then $waw$ contains a
$\frac{7}{3}$-power $x$, where $|x| \leq 5 \cdot 2^i$.
\end{lemma}

\begin{proof}
Suppose $a=1$ (the case $a=0$ follows similarly).
The proof is by induction on $i$.  For the base case we have
$i=0$ and $|w|=4$.  An easy computation suffices to verify that
for all $w$ with $|w|=4$, $w1w$ contains a $\frac{7}{3}$-power $x$,
where $|x| \leq 5$ as required.

Let us now assume that the lemma holds for all $i'$, where $0<i'<i$.
Since $w$ avoids $\frac{7}{3}$-powers, and since $|w| \geq 7$, by
Theorem~\ref{KS2} we can write $w=u\mu(w')v$, where
$u,v\in\{\epsilon,0,1,00,11\}$ and $w'\in\Sigma_2^*$ is
$\frac{7}{3}$-power-free.  By applying
a case analysis similar to that used in Cases (1)--(4) of the
proof of Theorem~\ref{h_eq_mu} below, we can eliminate all but
three cases: $(u,v)\in\{(\epsilon,\epsilon),(\epsilon,0),(0,\epsilon)\}$.

\begin{enumerate}
\renewcommand{\labelenumi}{Case \arabic{enumi}:}
\item
$(u,v)=(\epsilon,\epsilon)$.
In this case $w=\mu(w')$.  This is clearly not possible, since for
$i>0$, $|w| = 5 \cdot 2^i - 1$ is odd.

\item
$(u,v)=(\epsilon,0)$.
Then $w=\mu(w')0$ and $w1w=\mu(w')01\mu(w')0=\mu(w'0w')0$.  If
$|w| = 5 \cdot 2^i - 1$, we see that $|w'| = 5 \cdot 2^{i-1} - 1$.
Hence, if $i'=i-1$, we may apply the inductive assumption to
$w'0w'$.  We thus obtain that $w'0w'$ contains a $\frac{7}{3}$-power
$x'$, where $|x'| \leq 5 \cdot 2^{i-1}$, and so $w1w$ must contain a
$\frac{7}{3}$-power $x=\mu(x')$, where $|x| \leq 5 \cdot 2^i$.

\item
$(u,v)=(0,\epsilon)$.  This case is handled similarly to the previous
case, and we omit the details.
\end{enumerate}

By induction then, we have that $waw$ contains a
$\frac{7}{3}$-power $x$, where $|x| \leq 5 \cdot 2^i$.
\end{proof}

\begin{lemma}
\label{pattern3}
For $i,j\in\mathbb{Z}^+$, let $w$ and $s$ be $\frac{7}{3}$-power-free
words over $\Sigma_2$ such that $|w| = 2^{i+1} - 1$ or
$|w| = 3 \cdot 2^i - 1$, and $|s| = 2^{j+1} - 1$ or
$|s| = 3 \cdot 2^j - 1$.  Assume also that $|s| \geq |w|$.
Let $a$ be an element of $\Sigma_2$. Then $sawawas$ contains a
$\frac{7}{3}$-power.
\end{lemma}

\begin{proof}
Suppose $a=1$ (the case $a=0$ follows similarly).
The proof is by induction on $i$.  For the base case we have
$i=1$ and either $|w|=3$ or $|w|=5$.  An easy computation
suffices to verify that for all $w$ with $|w|=3$ or $|w|=5$,
and all $a,b\in\Sigma_2^2$, $a1w1w1b$ contains a $\frac{7}{3}$-power.

Let us now assume that the lemma holds for all $i'$, where $1<i'<i$.
Since $w$ avoids $\frac{7}{3}$-powers, and since $|w| \geq 7$, by
Theorem~\ref{KS2} we can write $w=u\mu(w')v$, where
$u,v\in\{\epsilon,0,1,00,11\}$ and $w'\in\Sigma_2^*$ is
$\frac{7}{3}$-power-free.  Similarly, we can write $s=u'\mu(s')v'$, where
$u',v'\in\{\epsilon,0,1,00,11\}$ and $s'\in\Sigma_2^*$ is
$\frac{7}{3}$-power-free.  By applying a case analysis similar
to that used in Cases (1)--(4) of the proof of Theorem~\ref{h_eq_mu}
below, we can eliminate all but three cases: $(u,v,u',v')\in\{
(\epsilon,\epsilon,\epsilon,\epsilon),(\epsilon,0,0,\epsilon),
(0,\epsilon,\epsilon,0)\}$.

\begin{enumerate}
\renewcommand{\labelenumi}{Case \arabic{enumi}:}
\item
$(u,v,u',v')=(\epsilon,\epsilon,\epsilon,\epsilon)$.
In this case $w=\mu(w')$.  This is clearly not possible,
since for $i>1$, both $|w| = 2^{i+1} - 1$ and $|w| = 3 \cdot 2^i - 1$
are odd.

\item
$(u,v,u',v')=(\epsilon,0,\epsilon,0)$.
Then $w=\mu(w')0$, $s=\mu(s')0$, and $$s1w1w1s =
\mu(s')01\mu(w')01\mu(w')01\mu(s')0 = \mu(s'0w'0w'0s')0.$$  If
$|w| = 2^{i+1} - 1$ or $|w| = 3 \cdot 2^i - 1$, we see that
$|w'| = 2^i - 1$ or $|w| = 3 \cdot 2^{i-1} - 1$.  Similarly, if
$|s| = 2^{j+1} - 1$ or $|s| = 3 \cdot 2^j - 1$, we see that
$|s'| = 2^j - 1$ or $|s| = 3 \cdot 2^{j-1} - 1$.  Hence,
if $i'=i-1$, we may apply the inductive assumption to
$s'0w'0w'0s'$.  We thus obtain that $s'0w'0w'0s'$ contains a
$\frac{7}{3}$-power $x'$, and so $s1w1w1s$ must contain a
$\frac{7}{3}$-power $x=\mu(x')$.

\item
$(u,v,u',v')=(0,\epsilon,0,\epsilon)$.  This case is handled similarly
to the previous case, and we omit the details.
\end{enumerate}

By induction then, we have that $sawawas$ contains a $\frac{7}{3}$-power.
\end{proof}

\begin{lemma}
\label{form}
Let $n$ be a positive integer.  Then $n$ can be written in the form
$2^i - 1$, $3 \cdot 2^i - 1$, $5 \cdot 2^i - 1$,
or $(7+2j)2^i - 1$ for some $i,j\in\mathbb{N}$.
\end{lemma}

\begin{proof}
If $n=1$ then $n=2^1 - 1$ as required.  Suppose then that $n>1$.
Then we may write $n-1=m2^i$, where $m$ is odd and $i\in\mathbb{N}$.
But for any odd positive integer $m$, either $m\in\{1,3,5\}$, or $m$
is of the form $7+2j$ for some $j\in\mathbb{N}$, and the result follows.
\end{proof}

\section{Main theorem}
Let $h:\Sigma^* \rightarrow \Sigma^*$ be a morphism.  We say
that $h$ is \emph{non-erasing} if, for all $a\in\Sigma$,
$h(a)\neq\epsilon$.  Let $E$ be the morphism defined by $E(0)=1$
and $E(1)=0$.  The following theorem is analogous to a result
regarding overlap-free words due to Berstel and S\'{e}\'{e}bold
\cite{BS93}.

\begin{theorem}
\label{h_eq_mu}
Let $h:\Sigma_2^* \rightarrow \Sigma_2^*$ be a non-erasing morphism.
If $h(01101001)$ is $\frac{7}{3}$-power-free, then there exists an integer
$k \geq 0$ such that either $h=\mu^k$ or $h=E\circ\mu^k$.
\end{theorem}

\begin{proof}
Let $h(0)=x$ and $h(1)=x'$ with $|x|,|x'| \geq 1$.  The proof is by
induction on $|x|+|x'|$.  If $|x|<7$ and $|x'|<7$, then a quick computation
suffices to verify that if $h(01101001)$ is $\frac{7}{3}$-power-free,
then either $h=\mu^k$ or $h=E\circ\mu^k$, where
$k\in\{0,1,2\}$.  Let us assume then, without loss of generality,
that $|x| \geq |x'|$ and $|x| \geq 7$.  The word $x$ must avoid
$\frac{7}{3}$-powers, and so, by Theorem~\ref{KS2},
we can write $x=u\mu(y)v$, where $u,v\in\{\epsilon,0,1,00,11\}$
and $y\in\Sigma_2^*$.  We will consider all 25 choices for $(u,v)$.

\begin{enumerate}
\renewcommand{\labelenumi}{Case \arabic{enumi}:}
\item
$(u,v) \in \{(0,00),(00,0),(00,00),(1,11),(11,1),(11,11)\}$.
Suppose $(u,v)=(0,00)$.  Then $h(00)=0\mu(y)000\mu(y)00$ contains
the cube 000, contrary to the assumptions of the theorem.
The argument for the other choices for $(u,v)$ follows similarly.

\item
$(u,v) \in \{(0,11),(00,1),(00,11),(1,00),(11,0),(11,00)\}$.
For any of these choices for $(u,v)$, $h(00)=u\mu(y)vu\mu(y)v$ contains
a subword of the form $abb\mu(y)$ or $\mu(y)bba$ for some $a,b\in\Sigma_2$,
where $a \neq b$.  Since $|x| \geq 7$, $|y| \geq 2$, and so by
Lemma~\ref{abb} we have that $h(00)$ contains a $\frac{7}{3}$-power,
contrary to the assumptions of the theorem.

\item
$(u,v) \in \{(\epsilon,0),(0,\epsilon),(\epsilon,1),(1,\epsilon)\}$.
Suppose $(u,v)=(0,\epsilon)$.  Then $h(00)=0\mu(y)0\mu(y)$.  We have
two subcases.

\begin{enumerate}
\renewcommand{\labelenumii}{Case \arabic{enumi}\alph{enumii}:}
\item
$\mu(y)$ begins with 01 or ends with 10.  Then by Lemma~\ref{abb}, $h(00)$
contains a $\frac{7}{3}$-power, contrary to the assumptions of the theorem.

\item
$\mu(y)$ begins with 10 and ends with 01.  Then $h(00)=0\mu(y')01010\mu(y'')$
contains the $\frac{5}{2}$-power 01010, contrary to the assumptions
of the theorem.
\end{enumerate}

The argument for the other choices for $(u,v)$ follows similarly.

\item
$(u,v) \in \{(\epsilon,00),(0,0),(00,\epsilon),
             (\epsilon,11),(1,1),(11,\epsilon)\}$.
Suppose $(u,v)=(00,\epsilon)$.  Then $h(00)=00\mu(y)00\mu(y)$.
The word $\mu(y)$ may not begin with a 0 as that would create the
cube 000.  We have then that $h(00)=00\mu(y)0010\mu(y')$ for some
$y'\in\Sigma_2^*$. By Lemma~\ref{abb}, $h(00)$ contains a
$\frac{7}{3}$-power, contrary to the assumptions of the theorem.
The argument for the other choices for $(u,v)$ follows similarly.

\item
$(u,v) \in \{(0,1),(1,0)\}$.
Suppose $(u,v)=(0,1)$.  By Lemma~\ref{form}, the following three
subcases suffice to cover all possibilities for $|y|$.

\begin{enumerate}
\renewcommand{\labelenumii}{Case \arabic{enumi}\alph{enumii}:}
\item
$|y| = (7+2j)2^i - 1$ for some $i,j\in\mathbb{N}$.
We have $h(00)=0\mu(y)10\mu(y)1=0\mu(y1y)1$.  By Lemma~\ref{pattern},
$y1y$ contains a $\frac{7}{3}$-power.  The word $h(00)$ must then
contain a $\frac{7}{3}$-power, contrary to the assumptions of the theorem.

\item
$|y| = 5 \cdot 2^i - 1$ for some $i\in\mathbb{N}$.  Again
we have $h(00)=0\mu(y)10\mu(y)1=0\mu(y1y)1$.  By Lemma~\ref{pattern2},
$y1y$ contains a $\frac{7}{3}$-power.  The word $h(00)$ must then
contain a $\frac{7}{3}$-power, contrary to the assumptions of the theorem.

\item
$|y| = 2^i - 1$ or $|y| = 3 \cdot 2^i - 1$ for some
$i\in\mathbb{N}$.  We have two subcases.

\begin{enumerate}
\renewcommand{\labelenumiii}{Case \arabic{enumi}\alph{enumii}.\roman{enumiii}:}
\item
$|x'| < 7$.
We have $h(0110)=0\mu(y)1x'x'0\mu(y)1$.  The only $x'\in\Sigma_2^*$
where $|x'| < 7$ and $1x'x'0$ does not contain a $\frac{7}{3}$-power
is $$x'\in\{10,0110,1001,011010,100110,101001\}.$$  However, each of these
words either begins or ends with 10, and so we have that $h(0110)$
contains a subword of the form $100\mu(y)$ or $\mu(y)110$.  Hence, by
Lemma~\ref{abb} we have that $h(0110)$ contains a $\frac{7}{3}$-power,
contrary to the assumptions of the theorem.

\item
$|x'| \geq 7$.  By Theorem~\ref{KS2}, we can write $x'=u'\mu(z)v'$, where
$u',v'\in\{\epsilon,0,1,00,11\}$ and $z\in\Sigma_2^*$ is
$\frac{7}{3}$-power-free.  Applying the preceding case analysis to
$x'$ allows us to eliminate all but three subcases.

\begin{enumerate}
\renewcommand{\labelenumiv}
{Case \arabic{enumi}\alph{enumii}.\roman{enumiii}.\Alph{enumiv}:}
\item
$(u',v')=(0,1)$.
We have $$h(0110)=0\mu(y)10\mu(z)10\mu(z)10\mu(y)1=0\mu(y1z1z1y)1.$$
Moreover, by the same reasoning used in Case 5a and Case 5b,
we have $|z| = 2^j - 1$ or $|z| = 3 \cdot 2^j - 1$ for some
$j\in\mathbb{N}$, and so by Lemma~\ref{pattern3}, $y1z1z1y$ contains a
$\frac{7}{3}$-power.  The word $h(0110)$ must then contain a
$\frac{7}{3}$-power, contrary to the assumptions of the theorem.

\item
$(u',v')=(1,0)$.
Then $h(01)=0\mu(y)11\mu(z)0$.  The word $\mu(z)$ may not
begin with a 1 as that would create the cube 111.  We have then that
$h(01)=0\mu(y)1101\mu(z')0$ for some $z'\in\Sigma_2^*$.  By
Lemma~\ref{abb}, $h(01)$ contains a $\frac{7}{3}$-power, contrary
to the assumptions of the theorem.

\item
$(u',v')=(\epsilon,\epsilon)$.
Then $h(01)=0\mu(y)1\mu(z)$.  We have two subcases.

\begin{itemize}
\item
$\mu(z)$ begins with 01.  Then $h(01)=0\mu(y)101\mu(z')$ for some
$z'\in\Sigma_2^*$.  The word $\mu(y)$ may not end in 10 as that
would create the $\frac{5}{2}$-power 10101.  Hence
$h(01)=0\mu(y')01101\mu(z')$ for some $y'\in\Sigma_2^*$.
If $\mu(z')$ begins with 10, then $h(01)$ contains the $\frac{7}{3}$-power
0110110.  If $\mu(z')$ begins with 01, then $h(01)$ contains the
$\frac{5}{2}$-power 10101.  Either situation contradicts the assumptions
of the theorem.

\item
$\mu(z)$ begins with 10.  Then $h(01)=0\mu(y)110\mu(z')$ for some
$z'\in\Sigma_2^*$.  By Lemma~\ref{abb}, $h(01)$ contains a
$\frac{7}{3}$-power, contrary to the assumptions of the theorem.
\end{itemize}
\end{enumerate}
\end{enumerate}
\end{enumerate}

The argument for the other choice for $(u,v)$ follows similarly.

\item
$(u,v) = (\epsilon,\epsilon)$.
In this case we have $x=\mu(y)$.
\end{enumerate}

All cases except $x=\mu(y)$ lead to a contradiction.  The same reasoning
applied to $x'$ gives $x'=\mu(y')$ for some $y'\in\Sigma_2^*$.
Let the morphism $h'$ be defined by $h'(0)=y$ and $h'(1)=y'$.  Then
$h=\mu \circ h'$, and by Theorem~\ref{KS1}, $h'(01101001)$ is
$\frac{7}{3}$-power-free.  Moreover, $|y|<|x|$ and $|y'|<|x'|$.
Also note that for the preceding case analysis it sufficed to consider the
following words only: $h(00)$, $h(01)$, $h(10)$, $h(11)$, $h(0110)$,
$h(1001)$, and $h(01101001)$.  However, 00, 01, 10, 11, 0110, and 1001 are all
subwords of 01101001.  Hence, the induction hypothesis can be applied, and
we have that either $h'=\mu^k$ or $h'=E \circ \mu^k$.
Since $E \circ \mu = \mu \circ E$, the result follows.
\end{proof}

We now establish the following corollary.

\begin{corollary}
\label{equiv}
Let $h:\Sigma_2^* \rightarrow \Sigma_2^*$ be a morphism such that
$h(01)\neq\epsilon$.  Then the following statements are equivalent.
\begin{enumerate}
\renewcommand{\labelenumi}{(\alph{enumi})}
\item The morphism $h$ is non-erasing, and $h(01101001)$ is
$\frac{7}{3}$-power-free.
\item There exists $k \geq 0$ such that $h=\mu^k$ or $h=E\circ\mu^k$.
\item The morphism $h$ maps any infinite $\frac{7}{3}$-power-free word
to an infinite $\frac{7}{3}$-power-free word.
\item There exists an infinite $\frac{7}{3}$-power-free word whose
image under $h$ is $\frac{7}{3}$-power-free.
\end{enumerate}
\end{corollary}

\begin{proof}\

(a) $\Longrightarrow$ (b) was proved in Theorem~\ref{h_eq_mu}.

(b) $\Longrightarrow$ (c) follows from Lemma~\ref{KS1} via K\"{o}nig's
Infinity Lemma.

(c) $\Longrightarrow$ (d):  We need only exhibit an infinite
$\frac{7}{3}$-power-free word:  the Thue-Morse word, $\mu^\omega(0)$,
is overlap-free and so is $\frac{7}{3}$-power-free.

(d) $\Longrightarrow$ (a):  Let $\mathbf{w}$ be an infinite
$\frac{7}{3}$-power-free word whose image under $h$ is
$\frac{7}{3}$-power-free.  By Theorem~\ref{subw}, $\mathbf{w}$ must
contain 01101001, and so $h(01101001)$ is $\frac{7}{3}$-power-free.

To see that $h$ is non-erasing, note that if $h(0)=\epsilon$, then
since $h(01)\neq\epsilon$, $h(1)\neq\epsilon$.  But then
$h(01101001)=h(1)^4$ is not $\frac{7}{3}$-power-free, contrary to
what we have just shown.  Similarly, $h(1)\neq\epsilon$, and so
$h$ is non-erasing.
\end{proof}

Let $h:\Sigma_2^* \rightarrow \Sigma_2^*$ be a morphism.  We say
that $h$ is the \emph{identity morphism} if $h(0)=0$ and $h(1)=1$.
The following corollary gives the main result.

\begin{corollary}
\label{main}
An infinite $\frac{7}{3}$-power-free binary word is a fixed point
of a non-identity morphism if and only if it is equal to the Thue-Morse
word, $\mu^\omega(0)$, or its complement, $\mu^\omega(1)$.
\end{corollary}

\begin{proof}
Let $h:\Sigma_2^* \rightarrow \Sigma_2^*$ be a non-identity morphism, and
let us assume that $h$ has a fixed point that avoids $\frac{7}{3}$-powers.
Then $h$ maps an infinite $\frac{7}{3}$-power-free word to an infinite
$\frac{7}{3}$-power-free word, and so, by Corollary~\ref{equiv}, $h$
is of the form $\mu^k$ or $E\circ\mu^k$ for some $k \geq 0$.  Since $h$
has a fixed point, it is not of the form $E\circ\mu^k$, and since $h$
is not the identity morphism, $h=\mu^k$ for some $k \geq 1$.  But
the only fixed points of $\mu^k$ are $\mu^\omega(0)$ and $\mu^\omega(1)$,
and the result follows.
\end{proof}

\section{The constant $\frac{7}{3}$ is best possible}
\label{sec4}
It remains to show that the constant $\frac{7}{3}$ given in
Corollary~\ref{main} is best possible; \emph{i.e.}, Corollary~\ref{main}
would fail to be true if $\frac{7}{3}$ were replaced by any larger
rational number.  To show this, it suffices to exhibit an infinite
binary word $\mathbf{w}$ that avoids $(\frac{7}{3}+\epsilon)$-powers
for all $\epsilon > 0$, such that $\mathbf{w}$ is the fixed point of a morphism
$h:\Sigma_2^* \rightarrow \Sigma_2^*$, where $h$ is not of the form
$\mu^k$ for any $k \geq 0$.

For rational $\alpha$, we say that a word $w$ \emph{avoids $\alpha^+$-powers}
if $w$ avoids $(\alpha+\epsilon)$-powers for all $\epsilon > 0$.

Let $h:\Sigma_2^* \rightarrow \Sigma_2^*$ be the morphism defined by
\begin{eqnarray*}
h(0) & = & 0110100110110010110 \\
h(1) & = & 1001011001001101001.
\end{eqnarray*}
Since $|h(0)|=|h(1)|=19$, $h$ is not of the form $\mu^k$ for any $k \geq 0$.
We will show that the fixed point $h^\omega(0)$ avoids 
$\frac{7}{3}^+$-powers by using a technique similar to that given by
Karhum\"{a}ki and Shallit \cite{KS03}.  We first state the following lemma,
which may be easily verified computationally.

\begin{lemma}
\label{incl}
\begin{enumerate}
\renewcommand{\labelenumi}{(\alph{enumi})}
\item Suppose $h(ab)=th(c)u$ for some letters $a,b,c\in\Sigma_2$ and words
$t,u\in\Sigma_2^*$.  Then this inclusion is trivial (that is, $t=\epsilon$
or $u=\epsilon$).

\item Suppose there exist letters $a,b,c\in\Sigma_2$ and words
$s,t,u,v\in\Sigma_2^*$ such that $h(a)=st$, $h(b)=uv$, and $h(c)=sv$.
Then either $a=c$ or $b=c$.
\end{enumerate}
\end{lemma}

\begin{theorem}
\label{plus}
The fixed point $h^\omega(0)$ avoids $\frac{7}{3}^+$-powers.
\end{theorem}

\begin{proof}
The proof is by contradiction.  Let $w\in\Sigma_2^*$ avoid
$\frac{7}{3}^+$-powers, and suppose that $h(w)$ contains a
$\frac{7}{3}^+$-power.  Then we may write $h(w)=xyyy'z$
for some $x,z\in\Sigma_2^*$ and $y,y'\in\Sigma_2^+$,
where $y'$ is a prefix of $y$, and $|y'|/|y| > \frac{1}{3}$.
Let us assume further that $w$ is a shortest such string,
so that $0 \leq |x|,|z| < 19$.  We will consider two cases.

Case 1: $|y| \leq 38$.  In this case we have $|w| \leq 6$.  Checking
all 20 words $w\in\Sigma_2^6$ that avoid $\frac{7}{3}^+$-powers,
we see that, contrary to our assumption, $h(w)$ avoids
$\frac{7}{3}^+$-powers in every case.

Case 2: $|y| > 38$.  Noting that if $h(w)$ contains a
$\frac{7}{3}^+$-power, it must contain a square, we may apply a
standard argument (see \cite{KS03} for an example) to show that
Lemma~\ref{incl} implies that $h(w)$ can be written in the following form:
$$h(w) = A_1A_2 \ldots A_jA_{j+1}A_{j+2} \ldots A_{2j}A_{2j+1}A_{2j+2} \ldots
A_{n-1}A_n'A_n'',$$ for some $j$, where
\begin{eqnarray*}
A_i & = & h(a_i) \quad\mbox{for}\quad i=1,2,\ldots,n \quad\mbox{and}\quad
a_i\in\Sigma_2 \\
A_n & = & A_n'A_n'' \\
y   & = & A_1A_2 \ldots A_j \\
    & = & A_{j+1}A_{j+2} \ldots A_{2j} \\
y'  & = & A_{2j+1}A_{2j+2} \ldots A_{n-1}A_n' \\
z   & = & A_n''.
\end{eqnarray*}

Since $y'$ is a prefix of $y$, and since $|y'|/|y| > \frac{1}{3}$,
$A_n'$ must be a prefix of $A_k$, where $k = \lfloor\frac{j}{3}\rfloor + 1$.
However, noting that for any $a\in\Sigma_2$, any prefix of $h(a)$ suffices to
uniquely determine $a$, we may conclude that $A_k = A_n$.  Hence, we may write
$$h(w) = A_1A_2 \ldots A_{k-1}A_k \ldots A_jA_{j+1}A_{j+2} \ldots
A_{j+k-1}A_{j+k} \ldots A_{2j}A_{2j+1}A_{2j+2} \ldots A_{n-1}A_n,$$ where
\begin{eqnarray*}
y   & = & A_1A_2 \ldots A_{k-1}A_k \ldots A_j \\
    & = & A_{j+1}A_{j+2} \ldots A_{j+k-1}A_{j+k} \ldots A_{2j} \\
y'z & = & A_{2j+1}A_{2j+2} \ldots A_{n-1}A_n \\
    & = & A_1A_2 \ldots A_{k-1}A_k.
\end{eqnarray*}

We thus have $$w = (a_1a_2 \ldots a_j)^2a_1a_2 \ldots a_k,$$ where
$k = \lfloor\frac{j}{3}\rfloor + 1$.  Hence, $w$ is a $\frac{7}{3}^+$-power,
contrary to our assumption.  The result now follows.
\end{proof}

Theorem~\ref{plus} thus implies that the constant $\frac{7}{3}$ given in
Corollary~\ref{main} is best possible.

\section*{Acknowledgements}
The author would like to thank Jeffrey Shallit for suggesting the problem,
as well as for several other suggestions, such as the example
$001001\mu^\omega(1)$ given in the introduction, and for
pointing out the applicability of the proof technique used
in Section~\ref{sec4}.


\begin{thebibliography}{99}
\bibitem{ACS98}
J.-P. Allouche, J. Currie, J. Shallit, ``Extremal infinite overlap-free
binary words'', \emph{Electron. J. Combin.} \textbf{5} (1998), \#R27.
\bibitem{AS03}
J.-P. Allouche, J. Shallit, \emph{Automatic Sequences: Theory,
Applications, Generalizations}, Cambridge University Press, 2003.
\bibitem{BS93}
J. Berstel, P. S\'{e}\'{e}bold, ``A characterization of overlap-free
morphisms'', \emph{Discrete Appl. Math.} \textbf{46} (1993), 275--281.
\bibitem{Dej72}
F. Dejean, ``Sur un th\'{e}or\`{e}me de Thue'', \emph{J. Comb. Theory
Ser. A.} \textbf{13} (1972), 90--99.
\bibitem{GH64}
W. Gottschalk, G. Hedlund, ``A characterization of the Morse minimal set'',
\emph{Proc. Amer. Math. Soc.} \textbf{15} (1964), 70--74.
\bibitem{KS03}
J. Karhum\"{a}ki, J. Shallit, ``Polynomial versus exponential growth
in repetition-free binary words'' (2003). Preprint available at
{\tt http://www.arxiv.org/abs/math.CO/0304095}.
\bibitem{MH44}
M. Morse, G. Hedlund, ``Unending chess, symbolic dynamics, and a problem
in semi-groups'', \emph{Duke Math. J.} \textbf{11} (1944), 1--7.
\bibitem{RS84}
A. Restivo, S. Salemi, ``Overlap free words on two symbols''. In M. Nivat,
D. Perrin, eds., \emph{Automata on Infinite Words}, Vol. 192 of
\emph{Lecture Notes in Computer Science}, pp. 198--206, Springer-Verlag, 1984.
\bibitem{See82}
P. S\'{e}\'{e}bold, ``Morphismes it\'{e}r\'{e}s, mot de Morse et mot de
Fibonacci'', \emph{C. R. Acad. Sc. Paris} \textbf{295} (1982), 439--441.
\bibitem{See84}
P. S\'{e}\'{e}bold, ``Overlap-free sequences''. In M. Nivat, D. Perrin,
eds., \emph{Automata on Infinite Words}, Vol. 192 of \emph{Lecture Notes
in Computer Science}, pp. 207--215, Springer-Verlag, 1984.
\bibitem{See85}
P. S\'{e}\'{e}bold, ``Sequences generated by infinitely iterated morphisms'',
\emph{Discrete Appl. Math.} \textbf{11} (1985) 255--264.
\bibitem{Shu00}
A.M. Shur, ``The structure of the set of cube-free $\mathbb{Z}$-words in a
two-letter alphabet'' (Russian), \emph{Izv. Ross. Akad. Nauk Ser. Mat.}
\textbf{64} (2000), 201--224.  English translation in \emph{Izv. Math.}
\textbf{64} (2000), 847--871.
\bibitem{Thu12}
A. Thue, ``\"{U}ber die gegenseitige Lage gleicher Teile gewisser
Zeichenreihen'', \emph{Vidensk. I. Math. Nat. Kl.} \textbf{1} (1912), 1--67.
\end{thebibliography}
\end{document}